\documentclass[journal]{IEEEtran}

\usepackage{generic}
\usepackage{cite}
\usepackage{amsmath}
\usepackage{amssymb}
\usepackage{algorithmic}
\usepackage{graphicx}
\usepackage{textcomp}
\usepackage{dsfont}
\usepackage{breqn} 
\usepackage{algorithm}
\usepackage{algorithmic}
\usepackage{nicefrac}       
\usepackage{microtype}      
\usepackage{verbatim}
\usepackage{flushend}
\allowdisplaybreaks
\usepackage{xcolor}
\usepackage{soul}

\providecommand{\sep}{\ensuremath{\rm{sep}}}
\newcommand{\pard}[2]{\frac{\partial #1}{\partial #2}}
\newtheorem{theorem}{Theorem}[section]

\newtheorem{lemma}[theorem]{Lemma}
\newtheorem{definition}{Definition}[section]

\usepackage{algorithm}
\usepackage{algorithmic}
\usepackage{bbm}

\newcommand{\real}{{\mathbb{R}}} 
\usepackage{verbatim}
\def\BibTeX{{\rm B\kern-.05em{\sc i\kern-.025em b}\kern-.08em
    T\kern-.1667em\lower.7ex\hbox{E}\kern-.125emX}}

\begin{document}
\title{Distributionally Robust Surrogate Optimal Control for High-Dimensional Systems}
\author{Aaron Kandel,  Saehong Park,  and Scott J. Moura
\thanks{Submitted for review on June 10th, 2021. This work was supported by LG Chem Ltd., and in part by a National Science Foundation Graduate Research Fellowship.
}
\thanks{Aaron Kandel is affiliated with the Department of Mechanical Engineering at the University of California, Berkeley, Berkeley, CA 94709 USA (e-mail: aaronkandel@berkeley.edu).}
\thanks{Saehong Park and Scott Moura are affiliated with 
the Department of Civil and Environmental Engineering at the University of California, Berkeley, Berkeley, CA 94709 USA (e-mail: \{sspark, smoura\}@berkeley.edu).}}

\maketitle

\begin{abstract}
This paper presents a novel methodology for tractably solving optimal control and offline reinforcement learning problems for high-dimensional systems. This work is motivated by the ongoing challenges of safety, computation, and optimality in high-dimensional optimal control.  We address these key questions with the following approach. First, we identify a sequence-modeling surrogate methodology which takes as input the initial state and a time series of control inputs, and outputs an approximation of the objective function and trajectories of constraint functions.  Importantly this approach entirely absorbs the individual state transition dynamics. The sole dependence on the initial state means we can apply dimensionality reduction to compress the model input while retaining most of its information.  Uncertainty in the surrogate objective will affect the result optimality.  Critically, however, uncertainty in the surrogate constraint functions will lead to infeasibility, i.e. unsafe actions. When considering offline reinforcement learning, the most significant modeling error will be encountered on out-of-distribution data. Therefore, we apply Wasserstein ambiguity sets to ``robustify'' our surrogate modeling approach subject to worst-case out-of-sample modeling error based on the distribution of test data residuals. We demonstrate the efficacy of this combined approach through a case study of safe optimal fast charging of a high-dimensional lithium-ion battery model at low temperatures.  

\end{abstract}


\section{Introduction}
\label{sec:introduction}
\IEEEPARstart{T}{his} paper presents a novel model-based data-driven method for robust optimal control and offline reinforcement learning of high-dimensional dynamical systems.

Optimal control faces unique challenges related to guaranteeing optimality and computational efficiency \cite{Kirk00}.  These challenges are generally exacerbated when the dynamical system in question is a large-scale system, a classification based on the cardinality of state variables $n$ ($x\in\mathbb{R}^n$) being high (i.e. $n>10^2$ or $10^3$). Learning based methods can also struggle to guarantee feasible solutions.

In this work, we introduce a simple algorithmic framework which utilizes (i) neural function approximation, (ii) dimensionality reduction, and (iii) distributionally robust optimization (DRO) to obtain computationally tractable optimal control for large-scale nonlinear optimal control problems. This contribution is important, considering that the majority of real-life dynamical systems (i.e. heat transfer, fluid dynamics, etc...) are inherently large-scale.  This is partially a result of their representation with partial differential equations (PDEs), which when solved numerically are frequently represented with numerous state variables \cite{Biegler00}. Often, model-order reduction is applied to generate a ``control-oriented'' dynamical model when the true underlying system is complex and large-scale \cite{Kerschen00}.  However, reductions can refute our ability to observe fundamental insights from our optimal control solution \cite{Hespanha00}.  Reductions can also compromise the capability of maximizing the performance of the control policy. 

Relevant literature presents a host of methods for large-scale optimal control.  Besides use of specialized and case-specific heuristics, these generally include (i) control vector parameterization (CVP), (ii) reinforcement learning (RL) and approximate dynamic programming (ADP), (iii) pseudospectral optimal control (POC), and (iv) variational calculus and Pontryagin methods (PM) \cite{Moura00, Canon00}. 

CVP is a powerful tool due to its simplicity (see e.g. \cite{Methekar00}).  In CVP, the control input is represented and manipulated in reduced form.  For instance, the control input can be defined using a zero-order hold over long timesteps, or as a polynomial whose coefficients we optimize.  The advantage of CVP is it reduces the number of decision variables in the optimization program.   For instance, CVP has been used to reduce the complexity of highly non-convex but relatively small-scale problems \cite{Rothenberger00}.  Nonetheless, for large-scale control CVP has been shown to yield useful results \cite{Methekar00, Schlegel00}.  CVP simplifies the problem, which compromises optimality.  Furthermore, CVP only addresses computational cost from the cardinality of the control input. Other sources of computational expense (i.e. simulation, numerical optimization) can still prohibit tractable solution of the control problem.  


ADP leverages function approximation to enable policy learning beyond the spatial/memory limitations of tabular DP methods \cite{Bertsekas00, Bertsekas01}.   The three biggest shortcomings of ADP relate to safety, optimality, and computation.  ADP and other model-free RL methods often require constraints to be encoded as auxiliary penalties to the objective/reward function \cite{Garcia00}.  Weighting these penalties requires tuning the objective function carefully.  More importantly, however, model-free and model-based RL algorithms must learn behavior through exploration.  For constrained problems, this can implicitly require \textit{violating} constraints throughout online learning \cite{Ray00}.  Moreover, RL can lose guarantees of converging to an optimal policy when the problem is complex (i.e. \textit{not} linear-quadratic).  Furthermore, for large-scale nonlinear problems, ADP and model-free RL methods can require a large number of iterations to converge to a usable control policy  \cite{Bertsekas01}. At a high level, many of these challenges are just as relevant for model-based RL methods. These challenges are exacerbated when learning policies from fixed, offline datasets. Recent research in offline reinforcement learning literature has provided modified algorithms that address these challenges questions while also proving to be amenable to large-scale control \cite{nair2021awac, kumar2020conservative}. In particular, offline RL methods address distributional shifts between the training data and data encountered from novel experience. For high-dimensional systems, these shifts become more likely, and can hamper optimality and feasibility.

Surrogate optimization models typically map decision variables to an approximation of the true objective function.  Historically, surrogate optimization has been popular in aerospace applications, where complex high-dimensional physics-based models form the basis for design and analysis \cite{Mack00, Queipo00}.  The surrogate functions are fit using samples from the original objective, which is typically expensive to evaluate.  The most popular approach is efficient global optimization (EGO).  EGO is an adaptive sampling regime which is guaranteed to yield a surrogate optimization model with bounded modeling error under certain conditions \cite{Jones00}.  EGO can work for simple control problems \cite{Marzat00}, however for large-scale problems the parameterization of the surrogate model and the required sampling depth can become intractable.  Surrogate models have also been used to approximate state-transition dynamics for control \cite{Chen00,Nagabandi00}.  This application underpins modern research activity on model-based reinforcement learning \cite{Kaiser00, Landolfi00}.  For high-dimensional systems, such models are ostensibly impractical again due to the expansive parameterizations which would be required to represent state-transition dynamics.  


Table 1 shows a brief summary of the previously discussed algorithms. Existing methods possess unique strengths in solving large-scale optimal control problems, but there is area for further development.  The objective of this paper is to present a general, data-driven algorithmic framework applicable to high-dimensional systems which addresses the critical, unanswered question of safety and feasibility. First, we define neural network surrogates which map a reduced state representation and a finite time series of control inputs to an approximation of the objective function. Instead of constraint penalties, we develop auxiliary surrogate models which predict time series of the constraint functions using the same reduced input data. Our method is then, by definition, a model-based RL approach. For optimal control problems with a short time horizon, we obtain approximate solutions by optimizing around the models a single time.  However, for optimal control problems on longer time horizons, we apply these surrogates within a receding horizon control framework.   Via a sequence-modeling method, we absorb the dynamics of the state transitions into the prediction of the surrogate models, eliminating modeling drift.

\begin{table}[t]
\centering
\caption{Algorithms for High-Dimensional Control. A * indicates the approach can be applied given a fixed, offline dataset with no model knowledge.}
\label{table}
\setlength{\tabcolsep}{3pt}
\begin{tabular}{|p{35pt}|p{175pt}|} 
\hline
\textbf{Algorithm}& 
\textbf{Challenges}\\
\hline
CVP&
optimality, computation, requires model knowledge\\
RL*& 
safety, optimality, computation \\
POC& 
requires model knowledge, proprietary software \\
PM& 
numerical instability, computation, requires model knowledge \\
\hline
\end{tabular}
\label{tab1}
\end{table}


By leveraging surrogate models, we introduce modeling error.  While objective uncertainty may affect optimality, uncertainty in the constraint functions can mean the difference between safe control and critically unsafe behavior. Therefore, this work accommodates uncertainty in the constraints via distributionally robust chance constraints (DRCC). These chance constraints encode distributions of modeling error computed from testing data.  We apply Wasserstein ambiguity sets to strengthen robustness by optimizing with respect to worst-case modeling error sourced from a family of distributions within some Wasserstein distance of the empirical distribution.   The Wasserstein measure is distinguished from other probabilistic distances (i.e. moment-based methods of $\phi$-divergence \cite{Jiang00}) in that it is symmetric between two distributions, makes no assumptions on the shape of the distributions, and importantly provides an ``out-of-sample'' safety guarantee \cite{Esfahani00}.  When used for DRCCs, we can probabilistically guarantee adherence to constraints even when our surrogate models experience distributional shifts relative to the training data.  

To evaluate the efficacy of the algorithm, we solve the safe-fast charging problem for a high-dimensional lithium-ion battery model at low temperatures.  Lithium-ion battery fast charging is currently an active research area in the energy systems and controls literature.  Significant challenges can arise in this problem from using reduced-order models \cite{Perez00}.  If we leverage full order electrochemical battery models, then we benefit from more granular electrochemical information to safe operate the cell farther towards the boundary of its safe operating conditions \cite{Kandel01}. This increases the performance of the resulting charge/discharge cycle, but requires that we strictly adhere to safety constraints.  Violation of some electrochemical constraints leads to rapid aging and potential catastrophic cell failure.  Consequently, the fast charging problem presents a relevant safety-critical challenge to our proposed algorithm.  Historically, fast charging has been explored with reduced order models due to the nonlinearity and computational complexity of simulating the full-order dynamics \cite{Rahn2012,Canova2015, Forman2011b}.  By demonstrating that our surrogate optimal control algorithm can yield interesting charge cycles based on the full-order electrochemical model in real time, we validate its use for large-scale nonlinear optimal control problems. 

The results in this paper comprise a significant extension of our previous work in \cite{Kandel00}.  These extensions include (i) a comprehensive novel case study using a full-order electrochemical battery model, including a computational comparison to control using a reduced order model, and (ii) the use of Wasserstein ambiguity sets instead of more limited $\phi$-divergence.  

\begin{figure*}[]
      \centering  
      \includegraphics[trim = 0mm 0mm 0mm 0mm, clip, width=\textwidth]{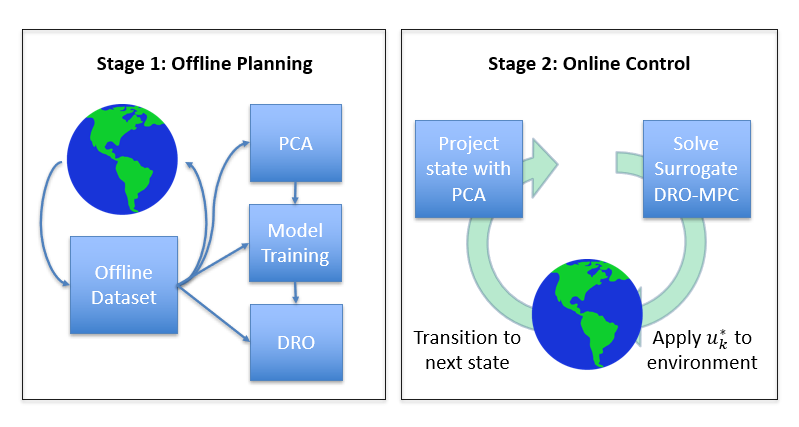}
      \caption{Block diagram detailing progression and flow of our proposed optimal control method. $u^*_k$ is the optimal open-loop control input obtained from MPC with the surrogate models. }
      \label{fig:blckdiag}
\end{figure*}

\section{Problem Formulation}
\subsection{Optimal Control Problem Formulation}

This paper considers the following optimal control problem statement, cast in discrete time: 
\begin{subequations}
\begin{align}
\min \quad & \sum_{k=0}^{N} J({x}_k,{u}_k) \label{eqn::ftocp1} \\
\text{subject to:} \quad
& {x}_{k+1} = f({x}_k,{u}_k) \label{eqn::ftocp2}\\
& g({x}_k,{u}_k) \leq 0 \label{eqn::ftocp3} \\
& h({x}_k,{u}_k) = 0 \label{eqn::ftocp4}\\
& {x}_0 = {x}(0) \label{eqn::ftocp5}
\end{align}
\end{subequations}
where $k$ is the current time and $N$ is the final time; ${x}_k \in \real^n$ is the state vector at time $k$; ${u}_k \in \real^p$ is the control input vector; $J({x}_k,{u}_k) : \real^n \times \real^p \rightarrow \real$ is the stage cost function at time $k$; $f({x}_k,{u}_k) : \real^n \times \real^p \rightarrow \real^n$ represents the system dynamics;  $g({x}_k,{u}_k) : \real^n \times \real^p \rightarrow \real^m$ represents inequality constraints; and $h({x}_k,{u}_k) : \real^n \times \real^p \rightarrow \real^\ell$ represents equality constraints. In this paper, we are particularly interested in problems where the cardinality of ${x}$ is high, i.e. $n > 10^2, 10^3$, or more. 

Our objective is to simplify the computation required to solve \eqref{eqn::ftocp1}-\eqref{eqn::ftocp5} when the model is high-dimensional. Figure \ref{fig:blckdiag} shows a block diagram of our method.  In each of the following sections, we discuss the components represented in this diagram.

\subsection{Offline Dataset}
Our method leverages a fixed, offline dataset composed of state trajectories matched with control input sequences. Typically, training data for surrogate optimization models is generated via a host of methods.  For instance, one popular method in the literature is Latin hypercube sampling (LHS) \cite{Queipo00}. In another method, EGO, sampling from the original objective function is organized and adaptive to the real-time evolution of modeling error \cite{Jones00}. In this paper, we train our surrogate models using data obtained from random, offline, parallelizable simulations of the original large-scale dynamical model. However, any dataset could be used to learn these surrogate models.  For example, such data could come from physical experiments, an existing suboptimal controller, etc...  In considering how such a dataset can be generated, the distributional shift problem becomes highly relevant.  We want to minimize the degree to which real-time control data will deviate from the distribution of training data. How this question is answered is highly dependent on the specific application.  Importantly, our framework is data-driven and does not require explicit model knowledge. This is differentiated from many existing methods (incl. CVP, psuedospectral optimal control).

\subsection{Model Formulation and Training}
Within the context of optimal control, surrogate models have been applied to represent state transition dynamics directly \cite{Chen00,Nagabandi00}. Direct approximation of state transition dynamics is not ideal for large-scale dynamical systems, where the large cardinality of state variables would require function approximators with intractable parameterizations. This paper proposes using a modified finite-time surrogate modeling approach which takes the following form:
\begin{subequations}
\begin{align}
\min \quad & \mathcal{J}({x}_0, U) \label{eqn:surr-ocp1} \\
\text{subject to:} \quad
& \mathcal{G}_i({x}_0, U) \leq 0 \: \forall \: i=1,\cdots,m \label{eqn:surr-ocp2}
\end{align}
\end{subequations}
The surrogate model $\mathcal{J}$ absorbs the state transition dynamics by mapping the initial state ${x}_0$ and time series of control inputs $U = [{u}(0), \cdots, {u}(N)]$ directly to an approximation of the objective function given in (\ref{eqn::ftocp1}). In set notation $\mathcal{J}(\cdot, \cdot) : \real^n \times \real^{p \times (N+1)} \rightarrow \real$. Likewise, the surrogate constraint functions $\mathcal{G}_i : \real^n \times \real^{p \times (N+1)} \rightarrow \real^{(N+1)}$ take the same inputs and predict as output a time series of the relevant constraint function values for each of $i = 1,...,m$ inequality constraints.   Importantly, the constraint surrogates only model the most relevant information in time series format. State variables that do not pertain to constraints in the optimization problem are disregarded by the surrogate models.  Furthermore, by outputting an entire time series, we avoid the possibility of modeling drift inherent to a surrogate which predicts individual state transitions across a single time step \cite{Kaiser00}.  

For a model predictive control application, the optimal control problem in \eqref{eqn:surr-ocp1}-\eqref{eqn:surr-ocp2} becomes:
\begin{subequations}
\begin{align}
\min \quad & \mathcal{J}({x}_k, U_{k:k+N}) \label{eqn:surr-opc-mpc1} \\
\text{subject to:} \quad &
\mathcal{G}_i({x}_k, U_{k:k+N}) \leq 0 \: \forall \: i=1,\cdots,m \label{eqn:surr-opc-mpc2}
\end{align}
\end{subequations}
At $k=0$, the initial state becomes the current state, and the control input time series $U_{k:k+N} = [{u}_k, \cdots, {u}(k+N)]$ starts at the current state and evolves over a horizon of $N$ time steps into the future. Note we are re-using $N$ here to indicate the control horizon length relative to the current time step, as opposed to the global time horizon length in \eqref{eqn:surr-ocp1}-\eqref{eqn:surr-ocp2}. After solving this reduced optimization program, we apply the first control input to the plant, simulating one step forward and then repeating the overall process.  

The most important transformation we make relates to reducing the state with dimensionality reduction techniques.  This paper specifically uses principal component analysis (PCA) to project the state onto a reduced basis. So in fact, the optimization program becomes:
\begin{subequations}
\begin{align}
\min \quad & \mathcal{J}(\tilde{x}_k, U_{k:k+N}) \label{eqn:surr-opc-mpc1} \\
\text{subject to:} \quad &
\mathcal{G}_i(\tilde{x}_k, U_{k:k+N}) \leq 0 \: \forall \: i=1,\cdots,m \label{eqn:surr-opc-mpc2}
\end{align}
\end{subequations}
where $\tilde{x}_k$ is a reduced representation of the dynamical state. Note the control is not included with state reduction, because its approximation could corrupt the input signal and negatively impact performance.

\subsubsection{Note: Facilitating Optimization}
This paper's approach requires we optimize around the neural network architecture. This architecture shares similar nonconvexity with the original expensive-to-evaluate objective function \cite{Jones00}. Past work has explored the use of convex neural architectures to facilitate this format of optimization \cite{Chen00}.  However, input-convex neural networks can compromise the universal function approximator properties of general neural networks \cite{Amos00}. 

The use of neural function approximation allows us to exploit analytic expressions for the function input-output gradient, as done in \cite{Chen00}.  For instance, for a single hidden layer neural network $f(x) = \sigma_{out}(W_2 \sigma_{hidden}(W_1 x + b_1) + b_2)$ where $\sigma_{out}(x) = x$, the Jacobian is given by:
\begin{equation}
    \text{Jac}(f(x))_{ij} = W_1(:,i)^T W_2(j,:) \sigma_{hidden}'(W_1 x + b_1 )
\end{equation}
Were we to solve the original optimal control problem with no surrogates, any gradients would be computed numerically via finite differences, which is highly inefficient. Numerical gradient calculations scale on the order of $\mathcal{O}[n^3]$ for a function $f: \real^n \rightarrow \real$, which would add significant computational complexity \cite{calafiore2014optimization}. By supplying the numerical optimization solver with analytic expressions for the input-output gradients of relevant surrogate models, we avoid expensive numerical gradient calculations.  Consequently, analytic gradients provide a fruitful opportunity to reduce computational complexity.

In this paper, we evaluate and compare two optimization schemes.  First, we use numerical optimization with specified analytical gradients.  We compare this approach to a sample-based random search.  Past work has shown for some applications that random search can provide high-performing results relative to more conventional optimization approaches \cite{Mania00}.  In this paper, we specifically apply a $(1+\lambda)$ evolutionary strategy algorithm to solve the receding horizon control problem. Section IV of this paper provides more details of this comparison. Overall, however, we were surprised by the ability of random search to outperform the gradient-based approach.

\section{Robustness to Modeling Errors}
Surrogate models are inherently imperfect.  Uncertainties are expected in approximations of both the objective and constraint functions and, if unaccounted for, these uncertainties can affect the optimality and feasibility of the final solution \cite{calafiore2014optimization}.

This paper addresses uncertainties in the constraint functions with a distributionally robust optimization (DRO) framework. We robustify our surrogate constraint models by optimizing with respect to worst-case realizations of modeling error characterized by the test data distribution of residuals. We obtain the worst-case realization through construction of a Wasserstein ambiguity set, which lends a probabilistic out-of-sample safety guarantee.  The following section details relevant mathematical preliminaries for this approach.

\subsection{Stochastic Optimization with Chance Constraints}
A chance constrained program includes probabilistic constraint statements, with random variables $\bf{R}$ with support $\Xi$.  Consider $x_k \in \mathbb{R}^{n}$ is the system state at time step $k$, $u(k) \in \mathbb{R}^{p}$ is the control input, $\textbf{R} \in \mathbb{R}^{m}$ is the random variable in question, and $g(x_k, u_k, \bf{R}): \mathbb{R}^n \times \mathbb{R}^p \times \mathbb{R}^m \rightarrow \mathbb{R}^m$ is the vector of inequality constraints.  The chance constraint is:
\begin{equation}\label{eqn:cc1}
    \hat{\mathbb{P}} \big{[}g(x_k, u_k, \textbf{R}) \leq 0\big{]} \geq 1 - \eta 
\end{equation}
where $\eta$ is our risk metric, or the probability of violating the constraint. The chance constraints discussed above depend on known distributions corresponding to each random variable.  For many applications, we approximate these distributions using data to create an empirical CDF.
In many data-driven applications, the true probability distribution $\mathbb{P}^*$ for the random variable $\bf{R}$ is unknown. Thus, our empirical distribution $\hat{\mathbb{P}}$ provides an approximation of $\mathbb{P}^*$ from data.  Borel's law of large numbers indicates that as the number of samples $\ell \rightarrow \infty$, $\hat{\mathbb{P}}\rightarrow \mathbb{P}^*$. This discrepancy characterizes distributional uncertainty in the random variable.  This can affect our solution if $\hat{\mathbb{P}}$ is inaccurate \cite{Nilim00}.  The literature presents several means by which we can accommodate this uncertainty.  In the following subsection, we discuss the application of the Wasserstein distance within this context.
\subsection{Wasserstein Ambiguity Sets}
An empirical distribution composed of samples will inevitably be characterized by some error or uncertainty.  In a qualitative sense, this uncertainty can be represented as the distribution lying some distance from the true distribution.
In statistics, there are several methods used to describe this type of distance.  These include $\phi$-divergence and the Wasserstein metric, the latter of which this paper applies for distributionally robust control.

\begin{definition}
Given two marginal probability distributions $\mathbb{P}_1$ and $\mathbb{P}_2$ lying within the set of feasible probability distributions $\mathcal{P}(\Xi)$, the Wasserstein distance between them is defined by
\begin{equation}
    \mathcal{W}(\mathbb{P}_1, \mathbb{P}_2) = \underset{\Pi} {\text{inf}} \bigg{\{} \int_{\Xi^2} ||\textbf{R}_1 - \textbf{R}_2 ||_a \Pi (d\textbf{R}_1, d\textbf{R}_2) \bigg{\}}
\end{equation}
where $\Pi$ is a joint distribution of the random variables $\bf{R}_1$ and $\bf{R}_2$, and $a$ denotes any norm in $\mathbb{R}^n$. 
\end{definition}

The Wasserstein distance allows us to replace the random variable with a ``worst-case'' realization sourced from a family of distributions within a specified Wasserstein distance of our empirical distribution. This family of distributions forms the Wasserstein ambiguity set.  For instance, let us define the ambiguity set as $\mathbb{B}_\epsilon$, a ball of probability distributions with radius $\epsilon$ centered around our empirical CDF $\hat{\mathbb{P}}$:
\begin{equation}\label{eqn:wass1}
    \mathbb{B}_\epsilon := \big{\{} \mathbb{P} \in \mathcal{P}(\Xi) \; | \; \mathcal{W}(\mathbb{P}, \hat{\mathbb{P}}) \leq \epsilon \big{\}}
\end{equation}
where $\epsilon$ is the Wasserstein ball radius. Now, we can formulate the robust counterpart of the chance constraint in (\ref{eqn:cc1}):
\begin{equation}\label{eqn:wass3}
    \underset{\mathbb{P} \in \mathbb{B}_\epsilon}{\text{inf}} \;  \mathbb{P} \big{[} g(x_k, u_k, \textbf{R}) \leq 0 \big{]} \geq 1 - \eta
\end{equation}
This equation provides the basis for the out of sample safety guarantee afforded by this DRO framework. Namely, we are probabilistically guaranteed to satisfy the chance constraint for any true probability distribution within $\epsilon$ distance of the empirical distribution.


Several expressions exist for the Wasserstein ball radius which, for a given confidence level $\beta$, is probabilistically guaranteed to contain the true distribution. We adopt the following formulation of $\epsilon$ from \cite{Zhao00} where $D$ is the diameter of the support of $\bf{R}$ composed of $\ell$ samples: 
\begin{equation}\label{eqn:wass2}
    \epsilon(\ell) = D \sqrt{\frac{2}{\ell} \log \bigg{(} \frac{1}{1-\beta} \bigg{)} }
\end{equation}
assuming the underlying data is independent and identically distributed (i.i.d.).
In \cite{Duan00}, this formula is replaced with the following version:
\begin{equation}\label{eqn:wass2}
    \epsilon(\ell) = C \sqrt{\frac{2}{\ell} \log \bigg{(} \frac{1}{1-\beta} \bigg{)} }
\end{equation}
where $C$ is obtained by solving the following scalar optimization program:
\begin{equation}\label{eqn:wass2}
    C \approx 2 \: \underset{\alpha > 0} {\text{inf}} \left\{ \frac{1}{2\alpha}\left( 1 + \ln \left(\frac{1}{N}\sum_{k=1}^N e^{\alpha ||\vartheta^{(k)}-\hat{\mu}||_1^2}\right)\right)\right\}^{\frac{1}{2}}
\end{equation}
where the right side bounds the value of $C$, and $\vartheta^{(k)}$ is a centered and normalized sample of the random variable which comprises our empirical distribution. This replacement is done to eliminate some unnecessary conservatism.

The exact constraint shown in (\ref{eqn:wass3}) is intractable, given that solving it requires solving an infinite dimensional nonconvex problem.  Most ongoing research in DRO focuses on deriving equivalent reformulations of (\ref{eqn:wass3}) which are more readily solved and implemented. 

What is particularly noteworthy about Wasserstein ambiguity sets is their inherent out-of-sample safety guarantee.  That is, the probabilistic safety guarantee covers, both, cases where we encounter experience and data which is not represented in the empirical distribution. This is principally due to the fact that the Wasserstein distance between two probability distributions bears no assumptions on the shape or support of each distribution.  We demonstrate this feature with the following comparison to $\phi$-divergence based reformulations of the constraint in (\ref{eqn:wass2}). 
If we were to utilize a $\phi$-divergence to reformulate (\ref{eqn:cc1}) as done in our previous work in \cite{Kandel00}:
\begin{equation}
    \mathbb{B}_\phi = \{\mathbb{P} \in \mathcal{P}(\Xi)\: | \: \phi(\mathbb{P}, \hat{\mathbb{P}}) \leq d \}
\end{equation}
where $d$ is a distance-like hyperparameter that must be tuned and chosen from intuition, then existing equivalent reformulations simply perturb the risk level \cite{Jiang00}.  However, perturbing the risk level provides much more limited out-of-sample guarantees because it limits the realization of the random variable to lie within a support that we have already observed. This finding is partially defined by the fact that the $\phi$-divergence between two probability distributions with different supports is infinite.  As a result, we adopt the Wasserstein distance metric for the remainder of this paper.  

\subsection{Equivalent Chance Constraint Reformulation}
In this paper, we adopt an equivalent reformulation of (\ref{eqn:wass3}) from \cite{Duan00}. This specific reformulation requires that the constraint function $g(x_k, u_k, \bf{R})$ is affine in $\bf{R}$. An in-depth discussion of this reformulation can be referenced in \cite{Duan00}.  Here, we \textit{restate} a brief overview of their methodology and derivation. 

We begin with samples of data $\{R^{(1)}, R^{(2)}, ..., R^{(\ell)} \}$ corresponding to random variable $\bf{R} \in \mathbb{R}^m$.  This sample comprises our empirical distribution $\hat{\mathbb{P}}$, and the data is drawn from the true underlying distribution $\mathbb{P}^*$.  First, we normalize the data samples to form a new random variable $\tilde{\vartheta}$ as follows:
\begin{equation}
    \vartheta^{(i)} = \Sigma^{-\frac{1}{2}}({R}^{(i)}-\mu)
\end{equation}
where $\Sigma$ is the sample variance of the data and $\mu$ is the sample mean. This standardization transforms the data samples such that its new mean is $0$, and its new variance is $I_{m\times m}$. Now, we define the support of this normalized distribution as:
\begin{equation}
    \Theta = \{\vartheta \in \mathbb{R}^m \ | \ -\sigma_{\max} \textbf{1}_m \leq \vartheta \leq \sigma_{\max} \textbf{1}_m\}
\end{equation}
Here, $\sigma_{\max} \in \mathbb{R}$ defines the support of the normalized random variable and $\bf{1}_m$ is a column vector of ones.  Now, let $\mathbb{Q}^*$ and $\hat{\mathbb{Q}}$ represent the true and empirical distributions of the normalized data $\vartheta$. We construct the ambiguity set ${\hat{\mathcal{Q}}}$ using the ``Wasserstein ball'' given by (\ref{eqn:wass1}), allowing us to transform the chance constraint in \eqref{eqn:cc1} to  
\begin{equation}
    \underset{\mathbb{Q} \in \hat{\mathcal{Q}}}{\text{sup}} \mathbb{Q}[ \vartheta \notin \mathcal{V}] \leq \eta
\end{equation}
which says the worst case probability that normalized random variable $\vartheta$ is outside set $\mathcal{V}$ is less than $\eta$, where the supremum is taken over all distributions $\mathcal{Q}$ in ambiguity set $\hat{\mathcal{{Q}}}$.
We wish to obtain the least conservative (i.e. tightest) set $\mathcal{V} \subseteq \mathbb{R}^m$ in order to define the desired Wasserstein uncertainty set $\mathcal{A} = \left\{ a \in \mathbb{R}^m \ | \ a = \Sigma^{\frac{1}{2}} v + \mu, \ v \in \mathcal{V} \right\}$ such that
\begin{equation}
    g(x_k, u_k, \bf{R}) \leq 0, \; \forall \; \bf{R} \in \mathcal{A}
\end{equation}
We restrict the overall shape of the set $\mathcal{V}$ to be a hypercube, which enables computational tractability:
\begin{equation}
    \mathcal{V}(\sigma) = \{ \vartheta \in \mathbb{R}^m | -\sigma \boldmath{1}_m < \vartheta < \sigma \boldmath{1}_m \}.
\end{equation}
Now, to compute this ambiguity set without introducing unnecessary conservatism, we need to find the minimum value of the hypercube side length $\sigma \in \mathbb{R}$.  The following optimization program details this problem:
\begin{align}
\underset{0\leq\sigma\leq \hat{\sigma}_{max}} {\text{min}}  &  \sigma \\ %
\text{subject to:} \quad &
\underset{\mathbb{Q} \in \hat{\mathcal{Q}}} {\text{sup}} \: \mathbb{Q}[\tilde{\vartheta} \notin \mathcal{V}(\sigma)]\leq \eta
\end{align}
Here, we select $\hat{\sigma}_{max}$ using \textit{a priori} information about the specific problem context.

The derivation in \cite{Duan00} provides a worst-case probability formulation, summarized by the following Lemma:
\begin{lemma}[Lemma 2 of \cite{Duan00}]
\begin{equation}
\begin{aligned} 
\underset{\mathbb{Q} \in \hat{\mathcal{Q}}} {\text{sup}} \mathbb{Q}[\tilde{\vartheta} \notin \mathcal{V}(\sigma)]=\\
\underset{\lambda \geq 0} {\text{inf}} \bigg{\{} \lambda \epsilon(\ell) + \frac{1}{\ell} \sum_{j=1}^\ell \left(1-\lambda \left(\sigma- ||\vartheta^{(j)}||_\infty \right)^+\right)^+\bigg{\}} \label{eqn:wass4}
\end{aligned}
\end{equation}
where $(x)^+=\max(x,0)$.
\end{lemma}

We defer to \cite{Duan00} for the proof of this finding. Their result entails that (\ref{eqn:wass4}) can be reformulated as
\begin{equation}\label{eqn:opt-reform}
     \underset{0 \leq \lambda,0\leq\sigma\leq \hat{\sigma}_{max}}  {\text{min}} \sigma \qquad \text{subject to:} \quad h(\sigma, \lambda) \leq \eta 
\end{equation}
where 
\begin{equation}
    h(\sigma, \lambda) = \lambda \epsilon(\ell) + \frac{1}{\ell} \sum_{j=1}^\ell \left( 1-\lambda(\sigma- ||\vartheta^{(j)}||_\infty)^+\right)^+
\end{equation}
The result of this optimization program is the value of $\sigma$, which is used to reformulate the chance constraints via convex approximation.  For a convex approximation of the constraint function in (\ref{eqn:wass3}), the hypercube $\mathcal{V}(\sigma)$ becomes the convex hull of its vertices.  If for example $m=1$ (i.e. the random variable is 1-dimensional), then $\mathcal{V}(\sigma)=(-\sigma, \sigma)$ -- an open interval. In general, this yields the ambiguity set $\mathcal{A} = \text{conv}(\{-r, r\})$ where $r=\Sigma^{\frac{1}{2}}\textbf{1}_m\sigma+\mu$ and $\text{conv}(\{\cdots\})$ represents the convex hull of points $\{\cdots\}$.  We can leverage this to complete the convex approximation of (\ref{eqn:wass3}) as a set of constraints of the form
\begin{equation}
    g(x_k, u_k) + r \leq 0
\end{equation}
which enumerate through the vertices of the robust hypercube.   For an $m$-dimensional constraint function, the exact form of the ambiguity set is $\mathcal{V} = \text{conv}(\{r^{(1)}, ..., r^{(2^m)}\})$.  The set of constraints are:
\begin{align}
    &g(x_k,u_k) + r^{(j)} \leq 0,  &\forall \ j=1,...,2^m
\end{align}
Algorithm 1 details the method used to compute the offset $\sigma$.
\begin{algorithm}[tb]
   \caption{Computation of $\sigma$}
   \label{alg:example}
\begin{algorithmic}
\item Initialize $\underline{\sigma}=0, \bar{\sigma}=\sigma_{max}$
\WHILE{$\bar{\sigma}-\underline{\sigma}>\text{tolerance}$}
   \item $\sigma = \frac{\bar{\sigma}+\underline{\sigma}}{2}$
   \item $[\lambda,h^*(\sigma,\lambda)]$ = minimize($\sigma$, $\lambda_{lb}$, $\lambda_{ub}$, $\epsilon$, $\theta$)
   \IF{$h^*(\sigma,\lambda)>\eta$}
   \item $\underline{\sigma}=\sigma$
   \ELSE
   \item $\bar{\sigma}=\sigma$
   \ENDIF
   \ENDWHILE
   \item $\sigma=\bar{\sigma}$
\end{algorithmic}
\end{algorithm}
In the next section, we detail exactly how we implement this robust optimization approach within the realm of learned optimal control, using a case study.

\subsection{Modeling Error as a Random Variable}
Notice that we have treated the residuals between the true and surrogate constraint functions, $\bf{R}$, as stochastic. In reality, the underlying process which generates these residuals can have deterministic structure, since they can be generated from deterministic models. That said, the training process yields an empirical set of residuals, for which an empirical probability distribution can be constructed. This stochastic modeling choice is convenient for chance constrained optimization, even if it neglects the underlying generative structure.

\section{Case Study}
Next we present a case study to validate and characterize the performance of the proposed algorithmic architecture.  Our case study is safe-fast charging of a lithium ion battery at low temperatures.  Lithium-ion battery fast charging is a highly relevant safety-critical application which possesses a rich and diverse history of research. It also presents a prototypical large-scale optimal control problem, in that complex electrochemical battery models are described with hundreds or even thousands of state variables.  While reduced-order equivalent circuit models address these dimensionality problems, the granular electrochemical information afforded by the full order models allows us to confidently take the battery closer to the safe operating envelope boundary.  This grants us the ability to exploit electrochemistry to improve charging performance \cite{Kandel00}.

Low temperatures complicate the fast charging problem problem, as they sensitize many of the complex electrochemical dynamics.  Specifically, the cell side-reaction overpotential constraint, which dictates the rate of lithium plating and cell degradation, can be much more readily violated at low temperatures \cite{Mohan16}.  Thus, the optimal control problem possesses many opportunities for constraint violation, which allows us to properly validate the efficacy of the proposed DRO framework.

Our case study is structured precisely as follows, where we solve a large-scale fast charging problem using the full-order Doyle-Fuller-Newman model (DFN) \cite{Thomas2002}. We also compare computation between the full order problem and one included in past work \cite{Kandel00} based on a moderately reduced single particle model.   We ensure comparison of our results with and without the added DRO framework, in order to validate its relative value and contributions to the safety of our algorithmic architecture.  


\subsection{Electrochemical Battery Model}
High fidelity battery modeling provides insights on performance, without requiring one to build and experimentally test the cell. The mathematical model formulated in this paper's appendix is the Doyle-Fuller-Newman battery model which comes from porous electrode theory, where Li-ions intercalate/deintercalate into porous spherical particles in the negative and positive electrodes. During charging, the Li-ions in the positive electrode deintercalate, dissolve into the electrolyte, and then migate and diffuse to the negative electrode by passing through the separator. Critically, this full-order electrochemical model reveals insights into the the mechanisms within the battery cell which allow us to take the battery farther towards the limit of its safe operating conditions.  By exploiting electrochemistry, we can calculate and apply faster, higher-performing charging cycles. 

While we relegate the model equations to this paper's appendix, we include some basic, useful information in this section in Table 1 for reference in discussing this paper's problem formulation and results. 
\begin{table}
\caption{Relevant Model Values}
\label{table}
\setlength{\tabcolsep}{3pt}
\begin{tabular}{|p{50pt}|p{100pt}|p{65pt}|}
\hline
\textbf{State Variable}& 
\textbf{Description}& 
\textbf{Units} \\
\hline
$SOC$& 
State of Charge& 
- \\
$\eta_S$& 
Side-Reaction Overpotential& 
Volts \\
$T$& 
Cell Temperature& 
 K \\
 $I$& Input Current & C-Rate\\
\hline
\end{tabular}
\label{tab1}
\end{table}


\subsection{Optimal Control Problem Statement}

For the DFN fast charging case study, we adopt the following optimal control problem statement within the framework of receding horizon control:
\begin{subequations}
\begin{align}
    & {\text{min}} \: \sum_{k=t}^{t+N} (SOC_k - SOC_{targ})^2 \\
    &\text{Subject to:} \\
    &\text{Dynamics}\\
    &\eta_s \geq 0 \\
    &T \leq T_{max} \\
    &0 \leq I \leq 2.5
    \end{align}
\end{subequations} 
The key constraints are that the side reaction overpotential stays positive, and the temperature does not exceed a maximum allowed threshold.  The overpotential constraint is the most critical barrier to prevent rapid aging and potential catastrophic failure of the cell.  If overpotential becomes negative, lithium metal begins to plate on the anode.  This phenomena reduces the capacity of the cell and leads directly to cell failure.  The temperature constraints provide indirect ways to avoid rapid aging, as the cell dynamics become more sensitive at temperature extremes.

We adapt this formulation using the distributionally robust surrogate modeling approach to yield:
\begin{subequations}
\begin{align}
    & {\text{min}} \: \mathcal{J}(x_{u,k}) \\
    &\text{subject to:} \\
    & \mathcal{G}_{\eta_s}(x_{u,k})  \geq r_{\eta_s} \\
    & \mathcal{G}_{T}(x_{u,k}) \leq T_{max}-r_{T}  \\
     &0 \leq I \leq 2.5
\end{align}
\end{subequations}
Since we are exploring fast charging at low temperatures, the temperature constraint is unlikely to be violated.  We omit this constraint, for simplicity, but it can be added back in practice.
\subsection{Results}
Table 2 details several important hyperparameters for this case study. We consider a nickel-manganese-cobalt battery cell. The initial electrochemical states correspond to equilibrium with a voltage of $V=3.25$ Volts.  The cell is at the same uniform temperature as the ambient temperature of $T_{amb}=281$ Kelvin.  We simulate 150 random charging trajectories to generate the requisite training data to fit the surrogate models. Each trajectory was either terminated if (1) the target SOC of 0.7 was reached, or (2) the episode end time of 55 minutes was reached. The maximum allowed C-rate for these simulations is 2.5C, where the C-rate for a lithium-ion battery is the parameter describing how much input current would be needed to charge the battery from empty to full in exactly 1 hour. A typical target SOC for electric vehicle applications is 0.8 or higher. Software implementations of the DFN model lose some numerical stability when applying high C-rates at higher SOCs.  To ensure we can continue utilizing a maximum C-rate of 2.5, we instead choose to set a slightly lower target SOC of 0.7 in our case study. Our algorithm can, however, be adapted to charge a battery cell to a higher SOC.

\begin{figure}[]
      \includegraphics[trim = 0mm 0mm 0mm 0mm, clip,scale=0.55]{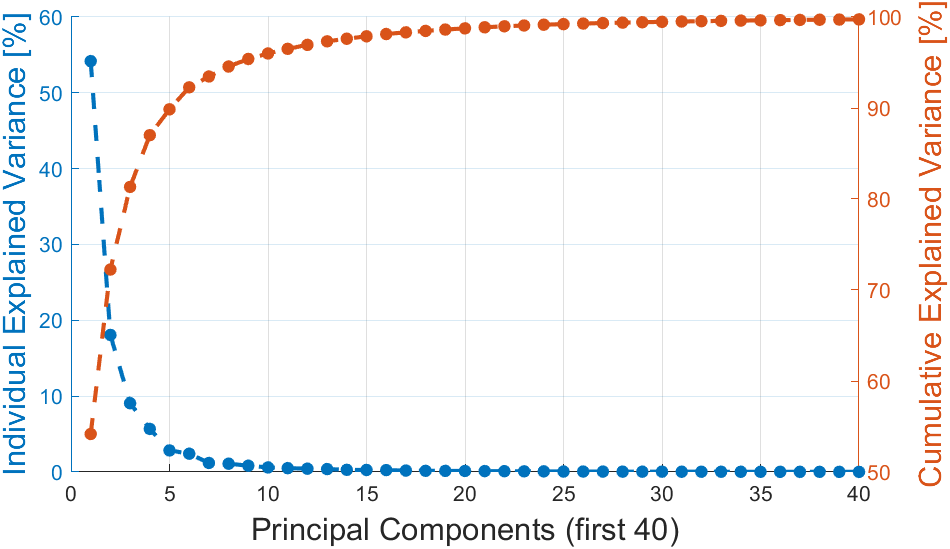}
      \caption{Individual and cumulative explained variance from principal component analysis of the electrochemical model state trajectories.}
      \label{fig:pca_var}
\end{figure}

Using principal component analysis on the state trajectories, we decide to project the state vector $x \in \real^{2687} \rightarrow \real^{40}$. This decision is motivated by the explained variance of the data, plotted in Fig. \ref{fig:pca_var}. Figure \ref{fig:pca_var} shows that the first 40 principal components of the state vector data explain 99.74\% of the variance in the dataset.

\begin{figure}[]
      \includegraphics[trim = 0mm 0mm 0mm 0mm, clip,scale=0.55]{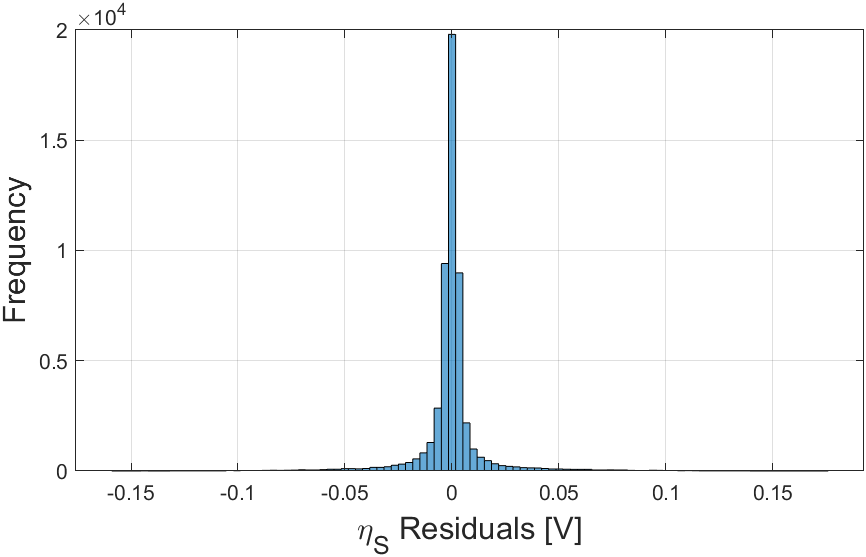}
      \caption{Histogram of test data residuals for $\mathcal{G}_{\eta_s}$.}
      \label{fig:pca_residuals}
\end{figure}

The surrogate models are feed-forward neural networks each with two hidden layers, each with 10 neurons and sigmoid activation functions.  The distribution of test data residuals for side reaction overpotential constraint function $\mathcal{G}_{\eta_s}$ are shown in Figure \ref{fig:pca_residuals}. This distribution is centered around zero with tight variance, although the tails of the distribution indicate that large residuals can occur with non-zero probability.  If unaccounted for in the control algorithm, violation of the overpotential constraint by, for example $0.14$ volts, would cause accelerated cell aging and could potentially sow the beginnings of a catastrophic failure.  Based on the testing data from model training (using an 80/20 split), the DRO offset computed using a Wasserstein ambiguity set is $r=0.0200066$. Given the specified chance constraint parameters, this offset is expected to yield desired safety characteristics.

\begin{table}
\caption{Relevant Hyperparameters}
\label{table}
\setlength{\tabcolsep}{3pt}
\begin{tabular}{|p{50pt}|p{100pt}|p{65pt}|}
\hline
\textbf{Parameter}& 
\textbf{Description}& 
\textbf{Value} \\
\hline
$\Delta$t& 
Timestep& 
15 seconds \\
$N$& 
Control Horizon& 
4 timesteps \\
$SOC_0$& 
Initial state-of-charge& 
0.0286 \\
$SOC_{targ}$& 
Target $SOC$& 
0.7 \\
$T_{amb}$& 
Ambient Temperature& 
281 Kelvin \\
$e$& 
Number of Training Episodes& 
100 \\
$T$& 
Length of Training Episode& 
3300 seconds \\
$I_{max}$& 
Maximum Charging Current& 
2.5 C \\
$\beta$& 
Ambiguity Set Confidence &  
0.9 \\
$\rho$& 
Chance Constraint Risk Metric& 
0.1 \\
\hline
\end{tabular}
\label{tab1}
\end{table}
\begin{figure*}[]
      \centering  
      \includegraphics[trim = 0mm 0mm 0mm 0mm, clip, width=\textwidth]{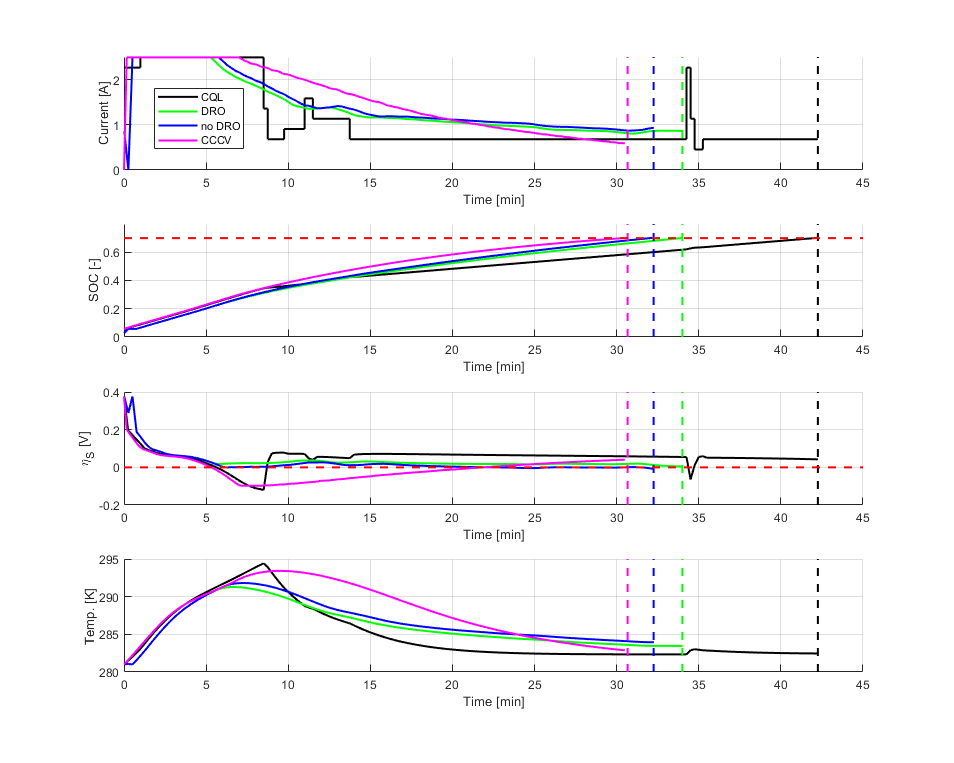}
      \caption{Optimal charging results for the DFN model using a nickel-manganese-cobalt (NMC) cell parameterization.  Here, the maximum allowed C-Rate is 2.5C and the target SOC is 0.7.  Charging is marked as complete at the vertical dotted lines for each respective trajectory.}
      \label{fig:finalresut}
\end{figure*}

We implemented our algorithm using a $(1+\lambda)$ evolutionary strategy for optimization, depending on 25000 mutants per iteration and 12 total iterations.  Cross-entropy random search also presents a useful alternative for numerical optimization \cite{BOTEV201335}. As a point of comparison, we implemented a numerical optimization scheme based on Matlab's \texttt{fmincon} solver, which we supplied with analytical gradient expressions for each function approximator.  The results from this implementation were inferior to a random search based optimization scheme.  The analytic gradients made \texttt{fmincon} nearly 70\% faster compared to using finite differences for gradient calculations. However, the average computation time per time step using \texttt{fmincon} was 9.1007 seconds whereas random search only required 2.0968 seconds per timestep on average.  We also find that the random search approach yields results of higher relative quality in terms of the overall charging time performance compared to the \texttt{fmincon} solver. The improved performance of random search, in terms of speed and solution equality, led us to use the random search method for our final results included in this paper. 

Our first benchmark is a hyper-aggressive constant current constant voltage (CCCV) charging protocol with $2.5$ C-rate maximum input current and $4.2$ Volts cutoff voltage.  A CCCV protocol charges the battery at the maximum allowed current until a cutoff voltage is reached.  From that point on, the battery is charged at a rate that keeps the voltage at the specified threshold. Typically, CCCV profiles correspond to thresholds given in the battery cell specifications document, which tend to limit the maximum allowed input current to around $1$C for most nickel-manganese-cobalt cells.  For the sake of consistency, we keep the maximum allowed current the same for each method. CCCV contextualizes the relative performance of the proposed method.

As a point of comparison, we also implement conservative Q-learning (CQL), a popular offline reinforcement learning algorithm that addresses distributional shift through penalties on out-of-distribution (OOD) actions \cite{kumar2020conservative}. The CQL network is a feed-forward network with two hidden layers each composed of 64 neurons, and ReLU activations.  The network input is the DFN state projected via the same PCA approach as our method. We discretize the input current into 11 bins between 0 and 2.5 C-rate. The network is trained in tandem with a target network iteratively with the same offline dataset used to learn the surrogate models of our approach. The reward function is given below, and is adopted with slight modification from recent work \cite{9668843} successfully applying actor-critic RL methods to lithium-ion battery fast charging:
\begin{equation}
    r = -I - 100(\textbf{1}_{\eta_S<0}|\eta_S|)
\end{equation}
A complementary OOD CQL loss is augmented to this reward function when training the networks \cite{kumar2020conservative}. CQL is a model-free method, meaning its sample efficiency isn't as high as our model-based approach.  In \cite{9668843},  model-free actor critic methods are shown to require on the order of $3e3$ episodes of learning to achieve high-performing charging results. Given in this case we are dealing with more than an order of magnitude reduction in available data, the fidelity of these CQL results is actually quite impressive.  CQL unfortunately does not provide certificates on safety and feasibility, which is reflected in the final charging profile as shown in Figure 4.  This highlights a comparative advantage of our model-based RL methodology, namely its out-of-sample safety guarantees.

\begin{figure}[]
    \centering 
      \includegraphics[scale=0.555]{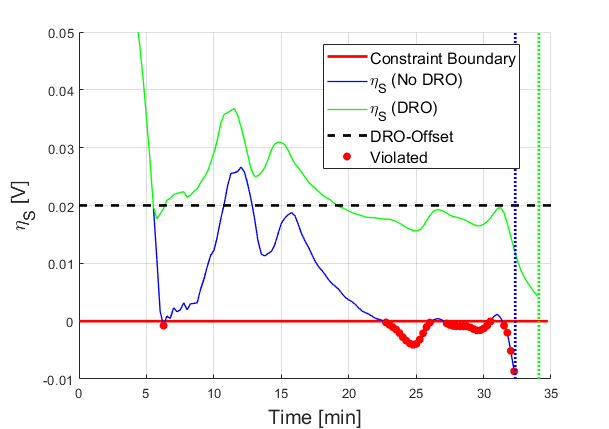}
      \caption{$\eta_s$ evolution statistics for each respective method (omitting CCCV and CQL). While the DRO version violates the conservative constraint offset, it still yields safe charging behavior relative to the nominal constraint boundary. Conversely, the non-robust version of our algorithm violates the nominal constraint boundary in 25.38\% of its timesteps.}
      \label{fig:etas}
\end{figure}

Figure \ref{fig:finalresut} shows the optimal fast charging results for versions of our algorithm with and without distributionally robust optimization.  Overall, the CCCV protocol charges in 30.6 minutes, the non-robust predictive controller in 32.35 minutes, the full distributionally robust controller in 34.1 minutes, and the CQL controller in 42 minutes.  The industry benchmark CCCV protocol yields a good performance with respect to charging time with a total time of 30.6 minutes.  However, it significantly violates the safety constraint by up to 0.12 Volts, and for extended periods of the overall experiment.  This would undoubtedly lead to significant degradation and potential failure of the cell. Figure \ref{fig:etas} shows constraint violation for each learning-based method.  Without the DRO architecture, the surrogate-based method provides a relatively high performing charging protocol which charges the battery cell in 32.35 minutes, only 5.7\% slower than the CCCV approach. It also demonstrates improved safety relative to the industry CCCV benchmark.  Specifically, the magnitude of the maximum constraint violation in the non-robust version of our algorithm is only 0.0082 Volts. With the added DRO framework based on Wasserstein ambiguity sets, we see that the charging protocol satisfies the constraint at every instance in time, while also providing a competitive 34.1 minute charging time.  These results illustrate the theoretical guarantees we expect from application of Wasserstein ambiguity sets. Relative to the non-robust version, the charging time with the DRO offset is only ~5.4\% slower, a tradeoff that may be worthwhile for the increased safety and mitigation of aging. CQL violates overpotential constraints and charges slowly in comparison, however we trained the CQL network with the exact same dataset as used by our method for consistency.  An offline dataset with (i) more trajectories, and (ii) trajectories that more frequently violate constraints would yield higher performing CQL results, however such results would not have any guarantees of adhering to constraints.

\subsection{Computational Effort Analysis}
Comparing the computational requirements of this algorithm to those of our preliminary version in \cite{Kandel00} reveals a host of meaningful insights.  In this paper, we are performing optimal control on the DFN model, which is characterized by 2687 state variables.  In the past exploratory work, we tested a more rudimentary version of our algorithm on the single particle model with electrolyte and thermal dynamics (SPMeT), a model with 208 state variables. The average computation time per iteration with the DFN is 2.0968 seconds, when the algorithm is executed on a Windows desktop workstation equipped with a 9th generation Intel i5 processor. In \cite{Kandel00}, the average time per iteration was 1.7803 seconds when run on the same machine.  Despite the more than 10-fold increase in the cardinality of the state vector of each model, the computational effort of the proposed algorithm only changes marginally by 17.81\%. This slight difference is likely due to the more complex neural network architecture and DRO framework which we employ in our updated analysis.

\subsection{Insights from Wasserstein DRO Algorithm}

One unique aspect of this work from preliminary results presented in \cite{Kandel00} is the application of Wasserstein ambiguity sets. Wasserstein ambiguity sets are differentiated from $\phi$-divergence based chance constraint reformulation by their robust out-of-sample safety guarantee. We see this difference by observing that Wasserstein ambiguity sets provide a slightly more conservative result that that shown in previous work.  This finding is clear from our DFN case study.  The DRO does prevent constraint violation entirely compared to the non-robust version which only attenuates its magnitude relative to CCCV.  For safety critical control applications, this added safety from the out of sample safety guarantee is valuable.

\begin{figure}[]
      \includegraphics[trim = 0mm 0mm 0mm 0mm, clip,scale=0.555]{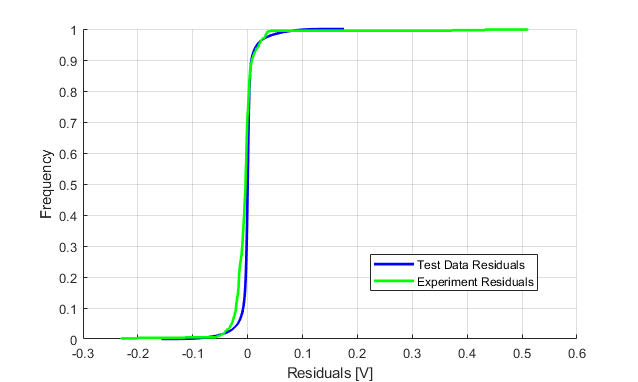}
      \caption{Comparison of cumulative distribution of $\mathcal{G}_{\eta_S}$ model residuals from test data and from the final optimal charging profile. These differences visualize the distributional shift problem that is a critical challenge in offline reinforcement learning.}  
      \label{fig:cdfcomp}
\end{figure}

To further demonstrate this added value, we refer to Figure \ref{fig:cdfcomp} which shows a comparison of the cumulative distribution of $\mathcal{G}_{\eta_S}$ model residuals from test data and from the state-action pairs in the final optimal charging profile. This plot highlights the distributional shift problem which is a significant open challenge in offline RL research. Consider that when limited to a static, offline dataset for model training, applying resulting control policies to a real, dynamical system creates the opportunity for the agent to encounter states that fall out of the distribution of its training data.  For high-dimensional, large-scale nonlinear dynamical systems, the probability of this occurring is significant. Thus, safety must be guaranteed with respect to such OOD experience.  Wasserstein ambiguity sets provide a strong means to satisfy this requirement, given their out-of-sample safety guarantee.  While the final experimental distribution does not represent the true underlying distribution of residuals, it does present a significant deviation from what we observe in our test data. Besides some slight differences in overall shape, the experimental residual distribution is more heavily skewed to higher magnitudes of modeling error.  Importantly in this case the maximum residual we observe is 0.5033 Volts, which is 2.908 times the magnitude of the largest residual represented in the test data set.  This difference is just one way of demonstrating how distributional errors can come into play once we set out to apply an optimal charging policy.


\section{Conclusion}

This paper presents a novel framework for optimal control of high-dimensional large scale dynamical systems.  The key challenges to numerical optimal control addressed by this paper include: (i) the ``curse of dimensionality'' incurred by high-dimensional systems, (ii) formulations that are not linear-quadratic, and (iii) ensuring safety/feasibility when constraint model errors occur.  

We identify surrogate models that learn from limited offline datasets, and which absorb state transition dynamics to reduce compounded modeling error.  Principal component analysis applied to the training data allows us to project the high-dimensional data onto a reduced basis.  This makes the modeling architecture conducive to fast identification and evaluation.  Finally, we integrate these models into a receding horizon control framework.  Critically, our strategy utilizes distributionally robust optimization to robustify the solution to errors in the constraint function surrogate models. the OOD safety guarantee of Wasserstein DRO directly addresses the open challenge of distributional shift for offline RL problems.  All combined, we demonstrate that the algorithmic approach yields tractable and robust control results for large-scale dynamical systems.


\bibliographystyle{./IEEEtran} 
\bibliography{./IEEEabrv,./generic-color}

\begin{thebibliography}{10}
\providecommand{\url}[1]{#1}
\csname url@rmstyle\endcsname
\providecommand{\newblock}{\relax}
\providecommand{\bibinfo}[2]{#2}
\providecommand\BIBentrySTDinterwordspacing{\spaceskip=0pt\relax}
\providecommand\BIBentryALTinterwordstretchfactor{4}
\providecommand\BIBentryALTinterwordspacing{\spaceskip=\fontdimen2\font plus
\BIBentryALTinterwordstretchfactor\fontdimen3\font minus
  \fontdimen4\font\relax}
\providecommand\BIBforeignlanguage[2]{{%
\expandafter\ifx\csname l@#1\endcsname\relax
\typeout{** WARNING: IEEEtran.bst: No hyphenation pattern has been}%
\typeout{** loaded for the language `#1'. Using the pattern for}%
\typeout{** the default language instead.}%
\else
\language=\csname l@#1\endcsname
\fi
#2}}

\bibitem{Kirk00}
D.~E. Kirk, \emph{Optimal Control Theory}.\hskip 1em plus 0.5em minus
  0.4em\relax Dover, 1970.

\bibitem{Biegler00}
L.~T. Biegler, O.~Ghattas, M.~Heinkenschoss, and B.~van Bloeman~Waanders,
  ``Large-scale pde-constrained optimization: an introduction,'' \emph{Lecture
  Notes in Computational Science and Engineering}, pp. 3--13, 2003.

\bibitem{Kerschen00}
G.~Kerschen, J.~Golinval, A.~Vakakis, and L.~Bergman, ``The method of proper
  orthogonal decomposition for dynamical characterization and order reduction
  of mechanical systems: An overview,'' \emph{Nonlinear Dynamics}, vol.~41,
  no.~1, 2005.

\bibitem{Hespanha00}
J.~Hespanha, \emph{Linear Systems Theory}.\hskip 1em plus 0.5em minus
  0.4em\relax Princeton University Press, 2009.

\bibitem{Moura00}
S.~Moura, N.~Chaturvedi, and M.~Krstic, ``Constraint management in li-ion
  batteries: A modified reference governor approach,'' in \emph{2013 American
  Control Conference}.\hskip 1em plus 0.5em minus 0.4em\relax Washington, DC
  USA: IFAC, IEEE, 2013.

\bibitem{Canon00}
M.~Canon, \emph{Theory of Optimal Control and Mathematical Programming}.\hskip
  1em plus 0.5em minus 0.4em\relax McGraw, 1970.

\bibitem{Methekar00}
R.~Methekar, V.~Ramadesigan, R.~Braatz, and V.~Subramanian, ``Optimum charging
  profile for lithium-ion batteries to maximize energy storage and
  utilization,'' \emph{Transactions of the Electrochemical Society}, vol.~25,
  no.~35, pp. 139--146, 2010.

\bibitem{Rothenberger00}
M.~J. Rothenberger, D.~J. Docimo, M.~Ghanaatpishe, and H.~K. Fathy, ``Genetic
  optimization and experimental validation of a test cycle that maximizes
  parameter identifiability for a li-ion equivalent-circuit battery model,''
  \emph{Journal of Energy Storage}, vol.~4, pp. 156--166, 2015.

\bibitem{Schlegel00}
M.~Schlegen, K.~Stockmann, T.~Binder, and W.~Marquardt, ``Dynamic optimization
  using adaptive control vector parameterization,'' \emph{Computers and
  Chemical Engineering}, vol.~29, no.~8, pp. 1731--1751, 2005.

\bibitem{Bertsekas00}
D.~P. Bertsekas, \emph{Dynamic Programming and Optimal Control}.\hskip 1em plus
  0.5em minus 0.4em\relax Athena Scientific Belmont, MA, 2017, vol.~1.

\bibitem{Bertsekas01}
D.~P. Bertsekas and J.~N. Tsitsiklis, \emph{Neuro-dynamic programming}.\hskip
  1em plus 0.5em minus 0.4em\relax Athena Scientific Belmont, MA, 1996, vol.~5.

\bibitem{Garcia00}
J.~Garcia and F.~Fernandes, ``A comprehensive survey on safe reinforcement
  learning,'' \emph{Journal of Machine Learning Research}, vol.~16, pp.
  1437--1480, 2016.

\bibitem{Ray00}
A.~Ray, J.~Achiam, and D.~Amodei, ``Benchmarking safe exploration in deep
  reinforcement learning,'' \emph{arXiv}, 2020.

\bibitem{nair2021awac}
A.~Nair, A.~Gupta, M.~Dalal, and S.~Levine, ``Awac: Accelerating online
  reinforcement learning with offline datasets,'' 2021.

\bibitem{kumar2020conservative}
A.~Kumar, A.~Zhou, G.~Tucker, and S.~Levine, ``Conservative q-learning for
  offline reinforcement learning,'' 2020.

\bibitem{Mack00}
Y.~Mack, T.~Goel, W.~Shyy, and R.~Haftka, ``Surrogate model-based optimization
  framework: A case study in aerospace design,'' \emph{Evolutionary Computation
  in Dynamic and Uncertain Systems}, 2007.

\bibitem{Queipo00}
N.~Queipo, R.~Haftka, W.~Shyy, T.~Goel, R.~Vaidyanathan, and P.~K. Tucker,
  ``Surrogate-based analysis and optimization,'' \emph{Progress in Aerospace
  Sciences}, vol.~41, pp. 1--28, 2005.

\bibitem{Jones00}
D.~R. Jones, M.~Schonlau, and W.~J. Welch, ``Efficient global optimization of
  expensive black-box functions,'' \emph{Journal of Global Optimization},
  vol.~13, no.~1, pp. 455--492, 1998.

\bibitem{Marzat00}
J.~Marzat and H.~Piet-Lahanier, ``Design of nonlinear mpc by kriging-based
  optimization,'' in \emph{16th IFAC Symposium on System Identification}.\hskip
  1em plus 0.5em minus 0.4em\relax Brussels, Belgium: The International
  Federation of Automatic Control, 2012, pp. 1490--1495.

\bibitem{Chen00}
Y.~Chen, Y.~Shi, and B.~Zhang, ``Optimal control via neural networks: A convex
  approach,'' in \emph{International Conference on Learning Representations
  (ICLR)}, New Orleans, LA USA, 2019.

\bibitem{Nagabandi00}
A.~Nagabandi, G.~Kahn, R.~S. Fearing, and S.~Levine, ``Neural network dynamics
  for model-based deep reinforcement learning with model-free fine-tuning,'' in
  \emph{International Conference on Robotics and Automation (ICRA)}, Brisbane,
  Australia, 2018.

\bibitem{Kaiser00}
L.~Kaiser, M.~Babaeizadeh, P.~Milos, B.~Osinski, R.~H. Campbell, K.~Czechowski,
  D.~Erhan, C.~Finn, P.~Kozakowski, S.~Levine, A.~Mohiuddin, R.~Sepassi,
  G.~Tucker, and H.~Michalewski, ``Model based reinforcement learning for
  atari,'' \emph{arXiv}, 2019.

\bibitem{Landolfi00}
N.~C. Landolfi, G.~Thomas, and T.~Ma, ``A model-based approach for
  sample-efficient multitask reinforcement learning,'' \emph{arXiv}, 2019.

\bibitem{Jiang00}
R.~Jiang and Y.~Guan, ``Data-driven chance constrained stochastic programs,''
  \emph{Mathematical Programming}, vol. 140, no.~6, pp. 291--327, 2016.

\bibitem{Esfahani00}
P.~Esfahani and D.~Kuhn, ``Data-driven distributionally robust optimization
  using the wasserstein metric: Performance guarantees and tractable
  reformulations,'' \emph{Mathematical Programming}, vol. 171, no. 1--2, pp.
  115--166, 2018.

\bibitem{Perez00}
\BIBentryALTinterwordspacing
H.~E. Perez, N.~Shahmohammadhamedani, and S.~Moura, ``{Enhanced Performance of
  Li-Ion Batteries via Modified Reference Governors and Electrochemical
  Models},'' \emph{IEEE/ASME Transactions on Mechatronics}, vol.~20, no.~4, pp.
  1511--1520, August 2015. [Online]. Available:
  \url{https://ieeexplore.ieee.org/document/7004876}
\BIBentrySTDinterwordspacing

\bibitem{Kandel01}
A.~Kandel and S.~Moura, ``Safe wasserstein constrained deep q-learning,''
  \emph{arXiv}, 2020.

\bibitem{Rahn2012}
C.~D. Rahn and C.-Y. Wang, \emph{{Battery Systems Engineering}}.\hskip 1em plus
  0.5em minus 0.4em\relax John Wiley {\&} Sons, 2012.

\bibitem{Canova2015}
M.~Canova, K.~Pan, and G.~Fan, ``{A Comparison of Model Order Reduction
  Techniques for Electrochemical Characterization of Lithium-Ion Batteries},''
  in \emph{54th IEEE Conference on Decision and Control}, Osaka, Japan, 2015.

\bibitem{Forman2011b}
\BIBentryALTinterwordspacing
J.~C. Forman, S.~Bashash, J.~L. Stein, and H.~K. Fathy, ``Reduction of an
  electrochemistry-based li-ion battery model via quasi-linearization and
  pad{\'{e}} approximation,'' \emph{Journal of the Electrochemical Society},
  vol. 158, no.~2, pp. A93--A101, 2011. [Online]. Available:
  \url{http://jes.ecsdl.org/content/158/2/A93.abstract}
\BIBentrySTDinterwordspacing

\bibitem{Kandel00}
A.~Kandel, S.~Park, H.~E. Perez, G.~Kim, Y.~Choi, H.~J. Ahn, W.~T. Joe, and
  S.~Moura, ``Distributionally robust surrogate optimal control for large-scale
  dynamical systems,'' in \emph{Proceedings of the 2020 American Control
  Conference (to appear)}.\hskip 1em plus 0.5em minus 0.4em\relax Denver, CO
  USA: IEEE, 2020.

\bibitem{Amos00}
B.~Amos, L.~Xu, and J.~Z. Kolter, ``Input convex neural networks,'' in
  \emph{International Conference on Machine Learning (ICML)}, Sydney,
  Australia, 2017.

\bibitem{calafiore2014optimization}
G.~C. Calafiore and L.~El~Ghaoui, \emph{Optimization models}.\hskip 1em plus
  0.5em minus 0.4em\relax Cambridge university press, 2014.

\bibitem{Mania00}
H.~Mania, A.~Guy, and B.~Recht, ``Simple random search provides a competitive
  approach to reinforcement learning,'' \emph{arXiv}, 2018.

\bibitem{Nilim00}
A.~Nilim and L.~E. Ghaoui, ``Robust control of markov decision processes with
  uncertain transition matrices,'' \emph{Operations Research}, vol.~53, no.~5,
  2005.

\bibitem{Zhao00}
C.~Zhao and Y.~Guan, ``Data-driven risk-averse stochastic optimization with
  wasserstein metric,'' \emph{Operations Research Letters}, vol.~46, no.~2, pp.
  262--267, 2018.

\bibitem{Duan00}
C.~Duan, W.~Fang, L.~Jiang, L.~Yao, and J.~Liu, ``Distributionally robust
  chance-constrained approximate ac-opf with wasserstein metric,'' \emph{IEEE
  Transactions on Power Systems}, vol.~33, no.~5, pp. 4924--4936, 2018.

\bibitem{Mohan16}
S.~{Mohan}, Y.~{Kim}, and A.~G. {Stefanopoulou}, ``Energy-conscious warm-up of
  li-ion cells from subzero temperatures,'' \emph{IEEE Transactions on
  Industrial Electronics}, vol.~63, no.~5, pp. 2954--2964, 2016.

\bibitem{Thomas2002}
\BIBentryALTinterwordspacing
K.~Thomas, J.~Newman, and R.~Darling, ``{Mathematical modeling of lithium
  batteries},'' \emph{Advances in lithium-ion batteries}, pp. 345--392, 2002.
  [Online]. Available:
  \url{http://www.springerlink.com/index/RXM87M4067U87J65.pdf}
\BIBentrySTDinterwordspacing

\bibitem{BOTEV201335}
\BIBentryALTinterwordspacing
Z.~I. Botev, D.~P. Kroese, R.~Y. Rubinstein, and P.~L’Ecuyer, ``Chapter 3 -
  the cross-entropy method for optimization,'' in \emph{Handbook of
  Statistics}, ser. Handbook of Statistics, C.~Rao and V.~Govindaraju,
  Eds.\hskip 1em plus 0.5em minus 0.4em\relax Elsevier, 2013, vol.~31, pp.
  35--59. [Online]. Available:
  \url{https://www.sciencedirect.com/science/article/pii/B9780444538598000035}
\BIBentrySTDinterwordspacing

\bibitem{9668843}
S.~Park, A.~Pozzi, M.~Whitmeyer, H.~Perez, A.~Kandel, G.~Kim, Y.~Choi, W.~T.
  Joe, D.~M. Raimondo, and S.~Moura, ``A deep reinforcement learning framework
  for fast charging of li-ion batteries,'' \emph{IEEE Transactions on
  Transportation Electrification}, pp. 1--1, 2022.

\end{thebibliography}

\section*{Appendix}
\subsection{Doyle-Fuller-Newman Electrochemical Battery Model}
We consider the Doyle-Fuller-Newman (DFN) model to predict the evolution of lithium concentration in the solid $c_{s}^{\pm}(x,r,t)$, lithium concentration in the electrolyte $c_{e}(x,t)$, solid electric potential $\phi_{s}^{\pm}(x,t)$, electrolyte electric potential $\phi_{e}(x,t)$, ionic current $i_{e}^{\pm}(x,t)$, molar ion fluxes $j_{n}^{\pm}(x,t)$, and battery temperature $T(t)$. The x-dimension runs across the negative electrode, separator, and positive electrode. At each x-coordinate value in the negative and positive electrodes, we consider a particle where spherical lithium intercalation occurs. The governing equations in time are given by
\begin{dgroup*}\eqnumsep=10em
\begin{dmath}
	\frac{\partial c_{s}^{\pm}}{\partial t}(x,r,t) = \frac{1}{r^{2}}\label{eqn::dyn1} \frac{\partial}{\partial r} \left[ D_{s}^{\pm} r^{2} \frac{\partial c_{s}^{\pm}}{\partial r}(x,r,t) \right], \label{eqn:dfn_cs} 
\end{dmath}
\begin{dmath}
	 \varepsilon_{e}^{j} \frac{\partial c_{e}^{j}}{\partial t}(x,t) = \frac{\partial}{\partial x} \left[D_{e}^{\text{eff}}(c_{e}^{j}) \frac{\partial c_{e}^{j}}{\partial x}(x,t) + \frac{1 - t_{c}^{0}}{F} i_{e}^{j}(x,t) \right], \ \label{eqn:dfn_ce} 
\end{dmath}
\begin{dmath}
	m c_{P} \frac{dT}{dt}(t) =  \frac{1}{R_{th}} \left[ T_{\textrm{amb}} - T(t) \right] + \dot{Q},  \label{eqn:dfn_T} 
\end{dmath}
\end{dgroup*}
for $j \in \{-,\sep,+\}$ and $\dot{Q}$ is the rate of heat transferred to the system \cite{Thomas2002}, defined as
\begin{align}
    \dot{Q} & = I(t)\left[U^{+}(t) - U^{-}(t) - V(t) \right] - \\ &I(t)T(t)\frac{\partial}{\partial T}[U^{+}(t) - U^{-}(t)], \label{eqn:dfn_Q_dot}
\end{align}
and differential equations in space and algebraic equations are given by
\begin{dgroup*}\eqnumsep=4em
\begin{dmath}
\sigma^{\text{eff},\pm}  \cdot \frac{\partial \phi_{s}^{\pm}}{\partial x}(x,t) = i_{e}^{\pm}(x,t) - I(t), \label{eqn:dfn_phis} 
\end{dmath}
\begin{dmath}
\kappa^{\text{eff}}(c_{e}) \cdot \frac{\partial \phi_{e}}{\partial x}(x,t) = -i_{e}^{\pm}(x,t) + \kappa^{\text{eff}}(c_{e}) \cdot \frac{2RT}{F}(1 - t_{c}^{0}) \times \left(1 + \frac{d \ln f_{c/a}}{d \ln c_{e}}(x,t) \right) \frac{\partial \ln c_{e}}{\partial x}(x,t), \label{eqn:dfn_phie} 
\end{dmath}
\begin{dmath}
\frac{\partial i_{e}^{\pm}}{\partial x}(x,t) = a^{\pm} F j_{n}^{\pm}(x,t), \label{eqn:dfn_ie}
\end{dmath}
\begin{dmath}
j_{n}^{\pm}(x,t) =\frac{1}{F} i_{0}^{\pm}(x,t) \left[e^{\frac{\alpha_{a}F}{RT} \eta^{\pm}(x,t)} - e^{-\frac{\alpha_{c}F}{RT} \eta^{\pm}(x,t)} \right], \label{eqn:dfn_bv}
\end{dmath}
\begin{dmath}
i_{0}^{\pm}(x,t) = k^{\pm}  \left[ c_{ss}^{\pm}(x,t) \right]^{\alpha_{c}} \left[c_{e}(x,t) \left(c_{s,\max}^{\pm} - c_{ss}^{\pm}(x,t)  \right) \right]^{\alpha_{a}}, \label{eqn:dfn_i0}
\end{dmath}
\begin{dmath}
\eta^{\pm}(x,t) = \phi_{s}^{\pm}(x,t) - \phi_{e}(x,t) - U^{\pm}(c_{ss}^{\pm}(x,t)) - F R_{f}^{\pm} j_{n}^{\pm}(x,t), \label{eqn:dfn_eta} 
\end{dmath}
\begin{dmath}
c_{ss}^{\pm}(x,t)= c_{s}^{\pm}(x,R_{s}^{\pm},t). \label{eqn:dfn_css}
\end{dmath}
\end{dgroup*}
\noindent where $D_{e}^{\text{eff}} = D_{e}(c_{e}) \cdot (\varepsilon_{e}^{j})^{\text{brug}}$, $\sigma^{\text{eff}} = \sigma \cdot (\varepsilon_{s}^{j} + \varepsilon_{f}^{j})^{\text{brug}}$, $\kappa^{\text{eff}} = \kappa(c_{e}) \cdot (\varepsilon_{e}^{j})^{\text{brug}}$ are the effective electrolyte diffusivity, effective solid conductivity, and effective electrolyte conductivity given by the Bruggeman relationship. The boundary conditions for solid-phase diffusion PDE (\ref{eqn:dfn_cs}) are
\begin{eqnarray}
	\frac{\partial c_{s}^{\pm}}{\partial r}(x,0,t) &=& 0, \label{eqn:dfn_cs-bc1}\\
	\frac{\partial c_{s}^{\pm}}{\partial r}(x,R_{s}^{\pm},t) &=& -\frac{1}{D_{s}^{\pm}} j_{n}^{\pm}(x,t). \label{eqn:dfn_cs-bc2}
\end{eqnarray}
The boundary conditions for the electrolyte-phase diffusion PDE (\ref{eqn:dfn_ce}) are given by
\begin{dgroup*}
\begin{dmath}
	\pard{c_{e}^-}{x}(0^{-},t) = {\pard{c_{e}^+}{x}(0^{+},t) = 0}, \label{eqn:dfn_ce-bc1}
\end{dmath}
\begin{dmath}
	 \varepsilon_{e}^{-} D_{e}(L^{-}) \pard{c_{e}^-}{x}(L^{-},t) = \varepsilon_{e}^{\textrm{sep}}D_{e}(0^{\textrm{sep}}) \pard{c_{e}^\sep}{x}(0^{\textrm{sep}},t), \ \ \ \ \label{eqn:dfn_ce-bc2}
\end{dmath}
\begin{dmath}
	\varepsilon_{e}^{\textrm{sep}}D_{e}(L^{\textrm{sep}}) \pard{c_{e}^\sep}{x}(L^{\textrm{sep}},t) = \varepsilon_{e}^{+} D_{e}(L^{+}) \pard{c_{e}^+}{x}(L^{+},t),  \label{eqn:dfn_ce-bc3}
\end{dmath}
\begin{dmath}
	c_{e}(L^{-},t) = c_{e}(0^{\textrm{sep}},t), \label{eqn:dfn_ce-bc4}
\end{dmath}
\begin{dmath}
	c_{e}(L^{\textrm{sep}},t) = c_{e}(L^+,t). \label{eqn:dfn_ce-bc5}
\end{dmath}
\end{dgroup*}
The boundary conditions for the electrolyte-phase potential ODE (\ref{eqn:dfn_phie}) are given by
\begin{eqnarray}
	\phi_{e}(0^{-},t) &=& 0,\label{eqn:dfn_phie-bc1} \\
	\phi_{e}(L^{-},t) &=& \phi_{e}(0^{\textrm{sep}},t), \label{eqn:dfn_phie-bc2}\\
	\phi_{e}(L^{\textrm{sep}},t)&=& \phi_{e}(L^+,t). \label{eqn:dfn_phie-bc3}
\end{eqnarray}
The boundary conditions for the ionic current ODE (\ref{eqn:dfn_ie}) are given by
\begin{equation}\label{eqn:dfn_ie-bcs}
	i_{e}^{-}(0^{-},t) = i_{e}^{+}(0^+,t) = 0,
\end{equation}
and also note that $i_{e}(x,t) = I(t)$ for $x \in [0^{\textrm{sep}}, L^{\textrm{sep}}]$. In addition, the parameters, $D_{s}^{\pm}, D_{e}, \kappa_{e}, k^{\pm}$ vary with temperature via the Arrhenius relationship:
\begin{align}
    \psi = \psi_{ref} \exp\left[ \frac{E_{\phi}}{R}\left( \frac{1}{T} - \frac{1}{T_{ref}} \right)\right]
\end{align}
where $\psi$ represents a temperature dependent parameter, $E_{\psi}$ is the activation energy and $\psi_{ref}$ is the reference parameter value at room temperature. The model input is the applied current density $I(t)$ [A/m$^{2}$], and the output is the voltage measured across the current collectors,
\begin{equation}\label{eqn:dfn_voltage}
	V(t) = \phi^{+}_{s}(0^{+},t) - \phi_{s}^{-}(0^{-},t).
\end{equation}

The level of charge in the cell is defined by the bulk state of charge (SOC) of the negative electrode, namely,
\begin{align}
    SOC^{-}(t) = \int_{0}^{L^{-}} \frac{\bar{c}_{s}^{-}(x,t)}{c_{s,max}^{-} (\theta^{-}_{100\%}-\theta^{-}_{0\%}) L^{-}} dx
\end{align}
where $\bar{c}_{s}^{-}$ represents the volume averaged of a particle in the solid phase defined as:
\begin{align}
    \bar{c}_{s}^{-}(x,t) = \frac{3}{(R_{s}^{-})^{3}}\int_{0}^{R_{s}^{-}} r^{2} c_{s}^{-}(r,t) dr
\end{align}

Lithium plating, which is the main battery degradation mechanism, is related to the side reaction overpotential $\eta_{s}$, defined as:
\begin{align} \label{eq:side_reaction_overpotential}
    \eta_{s}(x,t) = \phi_{s}^{-}(x,t) - \phi_{e}^{-}(x,t) - U_{sr} \geq 0.
\end{align}

To facilitate numerical optimal control, this model is discretized in space and time. There is a rich literature on discretization methods (see e.g. \cite{Rahn2012,Canova2015}). The discretization approached used for this paper involve finite difference, Pad\'{e} approximation \cite{Forman2011b}, and automatic differentiation methods.

\end{document}